\documentclass[reqno,11pt]{amsart}

\usepackage{amssymb,amsfonts, amsmath}
\usepackage{mathrsfs}
\usepackage{color}
\usepackage[normalem]{ulem} 

\newtheorem{thm}{Theorem}[section]
\newtheorem{cor}[thm]{Corollary}

\newtheorem{prop}[thm]{Proposition}
\newtheorem{rem}[thm]{Remark}

\numberwithin{equation}{section}
\newcommand{\Real}{\mathbb R}

\newcommand{\F}{\mathscr{F}}

\renewcommand{\P}{\mathsf{P}}

\newcommand{\E}{\mathsf{E}}

\DeclareMathOperator*{\essinf}{\mathrm{essinf}}
\DeclareMathOperator*{\sign}{\mathrm{sign}}
\newcommand{\Borel}{\mathscr{B}}

\begin{document}

\title[]{On a role of predictor in the filtering stability}%
\author{P. Chigansky}
\address{Department of Mathematics, The Weizmann Institute of Science, Rehovot 76100, Israel}
\email{pavel.chigansky@weizmann.ac.il}
\thanks{Research of the first author is supported by a grant from the Israel Science Foundation}

\author{R. Liptser}
\address{Department of Electrical Engineering Systems, Tel Aviv University, 69978 Tel Aviv, Israel}
\email{liptser@eng.tau.ac.il}



\begin{abstract}
When is a nonlinear filter stable with respect to its initial
condition? In spite of the recent progress, this question still
lacks a complete answer in general. Currently available results
indicate that stability of the filter depends on the signal
ergodic properties and the observation process regularity and may
fail if either of the ingredients is ignored. In this note we
address the question of stability in a particular weak sense and
show that the estimates of certain functions are always stable.
This is verified without dealing directly with the filtering
equation and turns to be inherited from certain one-step predictor
estimates.
\end{abstract}
\maketitle


\section{Introduction}
\label{sec-1}

Consider the filtering problem for a Markov chain
$(X,Y)=(X_n,Y_n)_{n\in\mathbb{Z}_+}$ with the signal $X$ and
observation $Y$. The signal process $X$ is a Markov chain itself
with the transition kernel $\Lambda(u,dx)$ and initial
distribution $\nu$. The observation process $Y$ has the transition
probability law
$$
\P(Y_n\in B|X_{n-1},Y_{n-1})=\int_B
\gamma(X_{n-1},y)\varphi(dy), \quad B\in\Borel(\Real),
$$
where $\gamma(u,y)$ is a density with respect to a $\sigma$-finite
measure $\varphi$ on $\Real$. We set $Y_0=0$, so that, a
priori information on the signal state at time $n=0$ is confined to
the signal distribution $\nu$. The random process $(X,Y)$ is
assumed to be defined on a complete probability space
$(\Omega,\F,\P)$. Let $(\F^Y_n)_{n\ge 0}$ be the filtration generated by $Y$:
$$
\F^Y_0=\{\varnothing,\Omega\}, \quad
\F^Y_n=\sigma\{Y_1,\ldots,Y_n\}.
$$
It is well known that the regular conditional distribution
$d\P(X_n\le x|\F^Y_n)=:\pi_n(dx)$ solves the
recursive Bayes formula, called the {\em nonlinear filter}:
\begin{equation}\label{filt}
\pi_n(dx) = \frac{\int_\Real\Lambda(u,dx) \gamma(u,Y_n)
\pi_{n-1}(du)} {\int_\Real\gamma(v,Y_n)\pi_{n-1}(dv)}, \quad n\ge
1,
\end{equation}
subject to $\pi_0(dx)=\nu(dx)$.
Clearly
$$
\pi_n(f):=\int_\Real f(x)\pi_n(dx)
$$
is a version of the
conditional expectation
$\E\big(f(X_n)|\F^Y_n\big)$ for any measurable
function $f=f(x)$, with $\E|f(X_n)|<\infty$.

Assume $\nu$ is unknown and the filter
\eqref{filt} is initialized with a probability distribution
$\bar{\nu}$, different from $\nu$ and denote the corresponding
solution by $\bar{\pi}=(\bar{\pi}_n)_{n\ge 0}$. Obviously, an
arbitrary choice of $\bar{\nu}$ may not be admissible: it makes
sense to choose $\bar{\nu}$ such that $\bar{\pi}_n(dx)$ preserves
the properties of a probability distribution, i.e.
$\int_B\bar{\pi}_n(dx)\ge 0$ for any measurable set
$B\in\Real$ and $\int_\Real \bar{\pi}_n(dx)=1$ for each
$n\ge 1$ $\P$-a.s. This would be the case if the right
hand side of \eqref{filt} does not lead to $0/0$ uncertainty with
a positive probability. As explained in the next section, the latter
is provided by the relation $ \nu\ll\bar{\nu}$, which is
assumed to be in force hereafter. In fact it plays an essential role in the proof
of main result.

The sequence $\bar{\pi}=(\bar{\pi}_n)_{n\ge 0}$ of random measures
generally differs from $\pi=(\pi_n)_{n\ge 0}$ and  the estimate $\bar{\pi}_n(f)$
of a particular function $f$ is said to
be stable if
\begin{equation}\label{wstab}
\E\big|\pi_n(f)-\bar{\pi}_n(f)\big|\xrightarrow[n\to\infty]{}0
\end{equation}
holds for any admissible pair $(\nu,\bar{\nu})$.

The verification of \eqref{wstab} in terms of $\Lambda(u,dx)$,
$\gamma(x,y)$, $\varphi(dy)$ is quite a nontrivial problem,
which is far from being completely understood in spite  of the
extensive research during the last decade.

For a bounded $f$,  \eqref{wstab} is closely related to ergodicity
of $\pi=(\pi_n)_{n\ge 0}$, viewed as a Markov process on the space
of probability measures. In the late 50's D. Blackwell, motivated
by the information theory problems, conjectured in \cite{Bl} that
$\pi$ has a unique invariant measure in the particular case of
ergodic Markov chain $X$ with a finite state space and noiseless
observations $Y_n=h(X_n)$, where $h$ is a fixed function. This
conjecture was found to be false by T. Kaijser, \cite{Kai}. In the
continuous time setting, H. Kunita addressed the same
question in \cite{K} for a filtering model with  general Feller-Markov
process $X$ and  observations
\begin{equation}\label{KKCL}
Y_t =\int_0^t h(X_s)ds + W_t,
\end{equation}
where the Wiener process $W=(W_t)_{t\ge 0}$ is independent
of $X$. According to \cite{K}, the filtering process
$\pi=(\pi_t)_{t\ge 0}$ inherits ergodic properties from $X$, if the tail $\sigma$-algebra of $X$ is $\P$-a.s. empty.
Unfortunately this assertion remains questionable due to a gap in
its proof (see \cite{BCL}).

\medskip
Notice that \eqref{wstab} for bounded $f$
also follows from
\begin{equation}\label{sstab}
\big\|\pi_n-\bar{\pi}_n\big\|_{\sf
tv}\xrightarrow[n\to\infty]{} 0, \quad \P-a.s.,
\end{equation}
where $\|\cdot\|_{\sf tv}$ is the total variation norm. Typically
this stronger type of stability holds when $X$ is an ergodic
Markov chain with the state space $\mathbb{S}\subseteq\Real$ (or
$\Real^d$, $d\ge 1$) and its transition probability kernel
$\Lambda(u,dx)$ is absolutely continuous with respect to a
$\sigma$-finite reference measure $\psi(dx)$,
$$
\Lambda(u,dx)=\lambda(u,x)\psi(dx),
$$
while the density $\lambda$ satisfies the so called {\em mixing}
condition:
\begin{equation}\label{mix}
0<\lambda_* \le \lambda(u,x)\le \lambda^*, \quad \forall \ x,u
\end{equation}
with a pair of positive constants $\lambda_*$ and $\lambda^*$.
Then (see \cite{AZ}, \cite{LGM}, \cite{DG}, \cite{CL}),
\begin{equation}\label{mmcc}
\varlimsup_{n\to\infty}\dfrac{1}{n}\log\big\|\pi_n-\bar{\pi}_n\big\|_{\sf
tv} \le -\dfrac{\lambda_*}{\lambda^*}, \quad \P\text{-}a.s.
\end{equation}
The condition \eqref{mix} was recently relaxed in \cite{CL}, where
\eqref{mmcc} was verified with $\lambda_*$ replaced by
$$
\lambda_\circ := \int_\mathbb{S} \essinf_{x\in
\mathbb{S}}\lambda(u,x)\mu(u)\psi(du),
$$
with $\mu(u)$ being the invariant density of the signal relative to
$\psi(du)$.

The mixing condition, including  its weaker form, implies
geometric ergodicity of the signal (see \cite{CL}). However, in
general the ergodicity (and even geometrical ergodicity) itself
does not imply stability of the filter (see
counterexamples in \cite{Kai}, \cite{DZ}, \cite{BCL}). If the
signal process $X$ is compactly supported, the density
$\lambda(u,x)$ usually corresponds to the Lebesgue measure
or purely atomic reference measure $\psi(dx)$. Signals with non
compact state space do not fit the mixing condition framework
since  an appropriate reference measure is hard to find and
sometimes it doesn't exist (as for the Kalman-Bucy filter).

In  non-compact or non-ergodic settings, the filtering
stability can be verified under additional structural assumptions
on $(X,Y)$. In this connection, we mention the
Kalman-Bucy filter being stable for
controllable and observable linear systems
 (see e.g. \cite{DZ}, \cite{OP}, \cite{LSII}, Sections 14.6 and 16.2).
Similarly, in the nonlinear case certain relations between $\lambda(x,u)$ and
$\gamma(x,y)$ provide  \eqref{sstab}  (see e.g. \cite{BO}, \cite{BO2}, \cite{A}, \cite{LO}, \cite{BCL}).

In summary, stability of the nonlinear filter stems from a delicate
interplay of the signal ergodic properties and the
observations ``quality''. If one of these ingredients is removed,
the other should be strengthened in order to keep the filter
stable. Notably all the available results verify
\eqref{wstab} via  \eqref{sstab} and, thus, require  restricting assumptions on
the signal structure.
Naturally, this
raises the following question: are there functions $f$ for which
\eqref{wstab} holds with ``minimal''
constraints on the signal model ?

In this note, we give examples of
functions for which this question has an
affirmative answer. It turns out that \eqref{wstab} holds if
$\nu\ll\bar{\nu}$ and the integral equation with respect to
$g$,
\begin{equation}\label{fff}
f(x)= \int_\Real g(y)\gamma(x,y)\varphi(dy),
\end{equation}
has  a bounded solution.
The proof of this fact relies on the martingale convergence
theorem rather than direct analysis of filtering equation \eqref{filt}.

The precise formulations and other generalizations with their  proofs are
given in Section \ref{sec-2}. Several nonstandard examples are
discussed in Section \ref{sec-3}.

\section{Preliminaries and the main result}\label{sec-2}

For notational convenience, we assume that the pair $(X,Y)$ is a
coordinate process defined on the canonical measurable space $(\Omega,\F)$ with
$\Omega=(\Real^\infty\times\Real^\infty)$ and
$\F=\Borel(\Real^\infty\times
\Real^\infty)$, where $\Borel$ stands for the Borel
$\sigma$-algebra. Let $\P$ be a probability measure on
$(\Omega,\F)$ such that $(X,Y)$ is   Markov
a process with the transition kernel
$
\gamma(u,y)\Lambda(u,dx)\varphi(dy)
$
and the initial distribution
$
\nu(dx)\delta_{\{0\}}(dy),
$
where $\delta_{\{0\}}(dy)$ is the point measure at zero.
Let $\bar{\P}$
be another probability measure on $(\Omega,\F)$ such that
$(X,Y)$ is  Markov process with the same transition
law  and the initial distribution
$\bar{\nu}(dx)\delta_{\{0\}}(dy)$. Hereafter,
$\E$ and $\bar{\E}$ denote expectations relative to $\P$ and $\bar{\P}$
respectively. By the Markov property of $(X,Y)$,
$$
\nu\ll\bar{\nu}\Rightarrow\P\ll\bar{\P}
\quad\text{and}\quad
\frac{d\P}{d\bar{\P}}(x,y)=\dfrac{d\nu}{d\bar{\nu}}(x_0),
\quad \bar{\P}\text{-a.s.}
$$
We assume that $\F^Y_0$ is completed with respect to
$\bar{\P}$. Denote $\F^Y_\infty =\bigvee_{n\ge 0}\F^Y_n$ and let $\P^Y$, $\bar{\P}^Y$
and $\P^Y_n$, $\bar{\P}^Y_n$ be the restrictions
of $\P$, $\bar{\P}$ on $\F^Y_\infty$ and
$\F^Y_n$ respectively. Obviously,
$\P\ll\bar{\P}\Rightarrow \P^Y\ll\bar{\P}^Y$ and $\P^Y_n\ll\bar{\P}^Y_n$ with the densities
$$
\frac{d\P^Y}{d\bar{\P}^Y}=
\bar{\E}\Big(\frac{d\nu}{d\bar{\nu}}(X_0)\big|\F^Y_\infty\Big)\quad
\text{and}\quad
\frac{d\P^Y_n}{d\bar{\P}^Y_n}=\bar{\E}
\Big(\dfrac{d\nu}{d\bar{\nu}}(X_0)\big|\F^Y_n\Big) :=
\varrho_n.
$$
Let $\bar{\pi}_n(dx)$ be the solution of \eqref{filt} subject to $\bar{\nu}$
considered on $\big(\Omega,\F,\bar{\P}\big)$, so that,
it is a version of the conditional distribution $\bar{\P}(X_n\le x|\F^Y_n)$.
Since $\P\ll\bar{\P}$, $\bar{\pi}_n$  satisfies \eqref{filt} on $\big(\Omega,\F,\P\big)$ as well.

In the sequel, we have to operate with $\varrho_n\pi_n(dx)$ as a random object
defined on $(\Omega,\F,\bar{\P})$. Since $\bar{\nu}\ll\nu$ is not assumed, $\pi_n$ cannot
be defined properly on $(\Omega,\F,\bar{\P})$ by applying the previous arguments.
However, the product $\varrho_n\pi_n$ is well defined on $(\Omega,\F,\bar{\P})$. Indeed, let
$\varGamma$ denote the set, where $\varrho_n\pi_n$ is well defined.
Notice that $\varGamma\in \F^Y_n$ and so,
$\bar{\P}(\varGamma)=\bar{\P}_n(\varGamma)$. Now,
by the Lebesgue decomposition of $\bar{\P}_n$ with respect to $\P_n$,
$$
\bar{\P}_n(\varGamma)=\int_{\varGamma\cap\{\varrho_n>0\}}
\varrho^{-1}_nd\P_n+\bar{\P}_n\big(\{\varrho_n=0\}\cap\varGamma
\big)
\ge \int_{\varGamma\cap\{\varrho_n>0\}} \varrho^{-1}_nd\P_n.
$$
Since both $\pi_n$ and $\varrho_n$ are defined $\P$-a.s.,
$\P_n(\varGamma)=1$ holds. Moreover, $\P_n(\varrho_n>0)=1$ since
$
\P_n(\varrho_n=0)=\int_{\{\varrho_n=0\}}\varrho_nd\bar{\P}_n=0.
$
Hence,
$$
\int_{\varGamma\cap\{\varrho_n>0\}} \varrho^{-1}_nd\P_n=
\int_\Omega \varrho^{-1}_nd\P_n
=\int_\Omega \varrho^{-1}_n\varrho_nd\bar{\P}_n=1,
$$
that is, $\bar{\P}_n(\varGamma)=1$.

\medskip
For $g:\Real\mapsto\Real$
with $\E|g(Y_n)|<\infty$ and
$\bar{\E}|g(Y_n)|<\infty$, let us define predicting
estimates:
$
\eta_{n|n-1}(g)=\E\big(g(Y_n)|\F^Y_{n-1}\big)
$
and
$
\bar{\eta}_{n|n-1}(g)=\bar{\E}\big(g(Y_n)|\F^Y_{n-1}\big).
$
We fix the following versions of these conditional expectations
$$
\begin{aligned}
\eta_{n|n-1}(g)&=\int_\Real \int_\Real
g(y)\gamma(x,y)\varphi(dy)\pi_{n-1}(dx)
\\
\bar{\eta}_{n|n-1}(g)&=\int_\Real \int_\Real
g(y)\gamma(x,y)\varphi(dy)\bar{\pi}_{n-1}(dx).
\end{aligned}
$$
Similarly to $\bar{\pi}_n$, the predictor $\bar{\eta}_{n|n-1}(g)$ is well defined
$\P$- and $\bar{\P}$-a.s. while only $\varrho_{n-1}\eta_{n|n-1}(g)$ makes sense
with respect to both measures.

\begin{thm}\label{lem-1}
Assume $\nu\ll\bar{\nu}$ and any of the following conditions:
\begin{enumerate}
\renewcommand{\theenumi}{{\bf \roman{enumi}}}
\item \label{i} $g$ is bounded{\rm ;}

\item \label{ii} $\dfrac{d\nu}{d\bar{\nu}}$ is bounded and the
family $(g(Y_n))_{n\ge 1}$ is $\bar{\P}$-uniformly
integrable{\rm ;}

\item \label{iii} for $p,q>1$, $\dfrac{1}{p}+\dfrac{1}{q}=1$,
$\bar{\E}\Big(\dfrac{d\nu}{d\bar{\nu}}\Big)^p<\infty$ and
the family $\big(|g(Y_n)|^q\big)_{n\ge 1}$ is
$\bar{\P}$-uniformly integrable.
\end{enumerate}
Then,
\begin{equation}\label{thm-claim}
\lim_{n\to\infty}\E\big|\eta_{n|n-1}(g)-\bar{\eta}_{n|n-1}(g)\big|=0.
\end{equation}
\end{thm}

\medskip
\noindent
{\bf Proof.} Suppose that $\alpha$ is $\F^Y_{n}$-measurable random variable defined
on $(\Omega,\F,\P)$ with $\E|\alpha|<\infty$. Then $\bar{\E}|\alpha|\varrho_{n}<\infty$ and
\begin{equation}\label{rrr}
\varrho_{n-1}\E\big(\alpha|\F^Y_{n-1}\big)=\bar{\E}
\big(\alpha\varrho_n|\F^Y_{n-1}\big), \quad
\bar{\P}  \text{-a.s.}
\end{equation}

{\bf (i)} For $\alpha:=g(Y_n)$, \eqref{rrr} reads:
$$
\varrho_{n-1}\eta_{n|n-1}(g)=
\bar{\E}\big(g(Y_n)\varrho_n|\F^Y_{n-1}\big),
\quad \bar{\P}\text{-a.s.}
$$
Therefore, (here $|g|\le C$ is assumed for definiteness)
\begin{equation}\label{2.44}
\begin{aligned}
\E\big|\eta_{n|n-1}(g)-\bar{\eta}_{n|n-1}(g)\big|
 &=
\bar{\E}\frac{d\P^Y_{n-1}}{d\bar{\P}^Y_{n-1}}
\big|\eta_{n|n-1}(g)-\bar{\eta}_{n|n-1}(g)\big|
\\
& =\bar{\E}\varrho_{n-1}
\big|\eta_{n|n-1}(g)-\bar{\eta}_{n|n-1}(g)\big|
\\
&= \bar{\E} \big|\varrho_{n-1}
\eta_{n|n-1}(g)-\varrho_{n-1}\bar{\eta}_{n|n-1}(g)\big|
\\
&= \bar{\E}
\big|\bar{\E}\big(g(Y_n)\varrho_n|\F^Y_{n-1}\big)
-\bar{\E}\big(g(Y_n)\varrho_{n-1}|\F^Y_{n-1}\big)\big|
\\
&=\bar{\E}
\Big|\bar{\E}\Big(g(Y_n)\big(\varrho_n-\varrho_{n-1}\big)\big|
\F^Y_{n-1}\Big)\Big| \le
 C\bar{\E}|\varrho_n-\varrho_{n-1}|.
\end{aligned}
\end{equation}
Since $(\varrho_n,\F^Y_n,\bar{\P})_{n\ge 1}$ is a
uniformly integrable martingale converging to
$
\varrho_\infty=\bar{\E}
\Big(\dfrac{d\nu}{d\bar{\nu}}\big|\F^Y_\infty\Big),
$
and $\bar{\E}|\varrho_n-\varrho_{n-1}|\le
\bar{\E}|\varrho_n-\varrho_\infty|+\bar{\E}|\varrho_{n-1} -\varrho_\infty|$,
the required result follows from $\lim_{n\to\infty}\bar{\E}|\varrho_n-\varrho_\infty|=0$ by
the Scheffe theorem.

\medskip
{\bf (ii)} Set $g^C=gI_{\{|g|\le C\}}$, then  by \eqref{i},
$$
\lim_{n\to\infty}\E\big|\eta_{n|n-1}(g^C)-\bar{\eta}_{n|n-1}(g^C)\big|=0,
\ \forall \ C>0.
$$
and it is left to show that
\begin{equation}\label{gC}
\begin{aligned}
&\lim_{C\to\infty}\varlimsup_{n\to\infty}\E\big|\eta_{n|n-1}(g-g^C)\big|=0
\\
&\lim_{C\to\infty}\varlimsup_{n\to\infty}\E\big|\bar{\eta}_{n|n-1}(g-g^C)\big|=0.
\end{aligned}
\end{equation}
Let for definiteness $\dfrac{d\nu}{d\bar{\nu}}\le K$ and thus
$\varrho_n\le K$, $\bar{\P}$-a.s. for all $n\ge 1$. Then
\begin{align*}
&\E\big|\eta_{n|n-1}(g-g^C)\big|\le
\E|g(Y_n)|I_{\{|g(Y_n)|>C\}} \le
K\bar{\E}|g(Y_n)|I{\{|g(Y_n)|>C\}}
\\
&\E\big|\bar{\eta}_{n|n-1}(g-g^C)\big|=\bar{\E}\varrho_{n-1}\big|
\bar{\eta}_{n|n-1}(g-g^C)\big| \le
K\bar{\E}|g(Y_n)|I{\{|g(Y_n)|>C\}},
\end{align*}
and \eqref{gC} holds by the uniform integrability assumption from \eqref{ii}.

\medskip
{\bf (iii)} By \eqref{2.44}, it suffices to show that
$\lim_{n\to\infty}\bar{\E}|g(Y_n)||\varrho_n-\varrho_{n-1}|=0.$ By the H\"older inequality we have
$$
\bar{\E}|g(Y_n)||\varrho_n-\varrho_{n-1}|\le
\Big(\bar{\E}|g(Y_n)|^q\Big)^{1/q}
\Big(\bar{\E}|\varrho_n-\varrho_{n-1}|^p\Big)^{1/p}.
$$
The $\bar{\P}$-uniform integrability of
$(|g(Y_n)|^q)_{n\ge 0}$ provides $\sup_{n\ge 0}\bar{\E}|g(Y_n)|^q<\infty$. Since
$$\lim_{n\to\infty}\bar{\E}|\varrho_n-\varrho_{n-1}|=0$$ it
is left to check that the family
$(|\varrho_n+\varrho_{n-1}|^p)_{n\ge 1}$ is
$\bar{\P}$-uniformly integrable. This holds by the
following upper bound
$$
\bar{\E}|\varrho_n+\varrho_{n-1}|^p\le
2^{p-1}(\bar{\E}\varrho^p_n+\bar{\E}\varrho^p_{n-1})\le
2^p \bar{\E}\Big(\dfrac{d\nu}{d\bar{\nu}}\Big)^p,\quad
p\ge 1
$$
where the Jensen inequality has been used.\hfill $\blacksquare$

\begin{cor}\label{theo-2}
Let $f$ be a measurable function and assume that there is a
function $g$ solving \eqref{fff} and satisfying the assumptions of
Theorem \ref{lem-1}. Then
$$
\lim_{n\to\infty}\E|\pi_n(f)-\bar{\pi}_n(f)|=0.
$$
\end{cor}

\medskip
\noindent
{\bf Proof.}
Since
$\pi_{n-1}(f)=\eta_{n|n-1}(g)\quad \text{and}\quad
\bar{\pi}_{n-1}(f)= \bar{\eta}_{n|n-1}(g)$, the claim is nothing
but \eqref{thm-claim}.\hfill $\blacksquare$

\section{Examples}\label{sec-3}

\subsection{Hidden Markov Chains} Let $X$ be a Markov chain taking values in a finite
alphabet $\mathbb{S}=\{a_1,...,a_d\}$ and the observation
$$
Y_n = \sum_{j=1}^d\xi_n(j)I_{\{X_{n-1}=a_j\}},
$$
where $\xi_n(j)$, $j=1,\ldots,d$, are independent entries of the random vectors
$\xi_n$, which form an i.i.d. sequence independent of $X$.

This variant of {\em Hidden Markov Model} is popular in various
applications (see e.g. \cite{EM}) and its stability analysis  has been
carried out by several authors (see e.g. \cite{AZ1}, \cite{LGM},
\cite{BCL}) mainly for ergodic chain $X$.
The nonlinear filter \eqref{filt} is finite dimensional, namely,
the conditional distribution $\pi_n(dx)$ is just the vector of
conditional probabilities
$\pi_n(i)=\P(X_n=a_i|\F^Y_n)$, $i=1,...,d$ and
$$
\|\pi_n-\bar{\pi}_n\|_{\sf
tv}=\sum_{i=1}^d|\pi_n(i)-\bar{\pi}_n(i)|.
$$

The following holds regardless of the ergodic properties of  $X$:
\begin{prop}\label{pro-4}
Assume
\begin{enumerate}
\renewcommand{\theenumi}{a\arabic{enumi}}
\item \label{a-1} all atoms of $\bar{\nu}$ are positive \item
\label{a-2} $\E\big|\xi_1(j)\big|^{i}<\infty$,
$i,j=1,\ldots,d$
\item \label{a-3} the $d\times d$ matrix $B$ with the entries $B_{ij}= \E \big(\xi_1(j)\big)^i$ is
nonsingular
\end{enumerate}
Then,
$$
\lim_{n\to\infty}\E \big\|\pi_n-\bar{\pi}_n\big\|_{\sf
tv}=0.
$$
\end{prop}

\medskip
\noindent
{\bf Proof.}
The condition \eqref{ii} of Theorem \ref{lem-1} is satisfied for
any $g_i(y)=y^i$, $i=1,\ldots,d$. Indeed, \eqref{a-1} and \eqref{a-2} imply
$\dfrac{d\nu}{d\bar{\nu}}\le\text{const.}$ and the uniform integrability of
$g_i(Y_n)$ for any $i$ since
$\bar{\E}\big|g_i(Y_n)\big|\le \sum_{j=1}^d\E\big|\xi_1(j)
\big|^i<\infty$.
Finally,
$$
\eta_{n|n-1}(g_i)=\E\big((Y_n)^i|\F^Y_{n-1}\big)
= \sum_{j=1}^d \pi_{n-1}(j)\E \big(\xi_1(j)\big)^i =
\sum_{j=1}^d \pi_{n-1}(j) B_{ij}.
$$
and, then, by Theorem \ref{lem-1},
$$
\E\big|\eta_{n|n-1}(g_i)-\bar{\eta}_{n|n-1}(g_i)\big| =
\E\Big|\sum_{j=1}^d
\big(\pi_{n-1}(j)-\bar{\pi}_{n-1}(j)\big)
B_{ij}\Big|\xrightarrow[n\to\infty]{}0.
$$
The latter and the nonsingularity of $B$ proves the claim.\hfill $\blacksquare$

\subsection{Observations with multiplicative white noise}

This example is borrowed from \cite{G-C}. The signal process is
defined by the linear recursive equation
$$
X_n = a X_{n-1} + \theta_n,
$$
where $|a|<1$ and $(\theta_n)_{n\ge 1}$ is $(0,b^2)$-Gaussian
white noise independent of $X_0$, that is, the signal process is
ergodic. The distribution function $\nu$ has density $q(x)$
relative to $dx$ from the {\em Serial Gaussian} (SG) family:
$$
q(x)=\bigg(\sum_{i=0}^\infty \alpha_i\frac{x^{2i}}{\sigma^{2i}C_{2i}}\bigg)\frac{1}{\sigma\sqrt{2\pi}}\exp\bigg(-\frac{x^2}{2\sigma^2}\bigg),
$$
where $\sigma$ is the scaling parameter, $\alpha_i$'s are nonnegative weight coefficients, $\sum_{i\ge 0}\alpha_i=1$ and
$C_{2i}$ are the normalizing constants.
The observation sequence is given by
$$
Y_n = X_{n-1}\xi_n,
$$
where $\xi_n$ is a sequence of i.i.d. random variables. The
distribution function of $\xi_1$ is assumed to have the following density
relative to $dx$:
\begin{equation}\label{ppp}
p(x)=\frac{\rho}{|x|^3}\exp\left(-\frac{\rho}{x^2}\right), \quad
p(0)=0,
\end{equation}
where $\rho$ is a positive constant. This filtering model is motivated by  financial applications when
$|X|$ is interpreted as the stochastic volatility parameter of an asset price.

As proved in \cite{G-C}, the filter \eqref{filt} admits
a finite dimensional realization provided that $\alpha_j\equiv 0$,
$j> N$ for some integer $N\ge 1$, namely for any time $n\ge 1$ the
filtering distribution $\pi_n(dx)$ has a density of SG type with
the scaling parameter $\sigma_n$ and the weights $a_{in}$, which
are propagated by a finite (growing with $n$) set of recursive
equations driven by the observations. Thus, the evolution of
$\pi_n(dx)$ is completely determined via $\sigma_n$ and
$\alpha_{in}$. Some stability analysis for the sequence
$\big(\sigma_n, (\alpha_{in})_{i\ge 1}\big)_{n\ge 1}$ has been
done in \cite{G-C2}.

Assume that the density $q(x)$ of $\nu$ belongs to the SG family, but its parameters are unknown.
If the filter is started from the Gaussian density with zero mean and variance ${\bar \sigma}^2$,
the filtering equation remains finite dimensional and the density
$$
\frac{d\nu}{d\bar{\nu}}(x)=\frac{q(x)}{\bar{q}(x)}=\bigg(\sum_{i=0}^N \alpha_i\dfrac{x^{2i}}{\sigma^{2i}C_{2i}}\bigg)
\frac{\bar \sigma}{\sigma} \exp\bigg(-\frac{x^2}{2}\Big(\frac{1}{\sigma^2}-\frac{1}{{\bar \sigma}^2}\Big)\bigg)
$$
is bounded, if $\bar\sigma>\sigma$.

In terms of the setting under consideration
$\gamma(x,y)=\dfrac{1}{|x|}p(y/x)$, where $p(\cdot)$ is defined in
\eqref{ppp} and $\varphi(dy)=dy$. For $f(x):=|x|$
$$
\eta_{n|n-1}(f)=\E\big(f(Y_n)|\F^Y_{n-1}\big)=\pi_{n-1}(f)\E
|\xi_1|,
$$
where $\E |\xi_1|>0$ and hence $ g(y)=|y|/\E |\xi_1| $
solves \eqref{fff}. Finally,  $\big(g(Y_n)\big)_{n\ge 1}$ is
$\bar{\P}$-uniformly integrable family since
$$
\begin{aligned}
\bar{\E}|X_{n-1}\xi_n|I_{\{|X_{n-1}\xi_n|>C\}}\le
\frac{\bar{\E}
|X_{n-1}\xi_n|^{1+\varepsilon}}{C^\varepsilon} =
\frac{\bar{\E}
|X_{n-1}|^{1+\varepsilon}\E|\xi_1|^{1+\varepsilon}}{C^\varepsilon},
\end{aligned}
$$
where $\E|\xi_1|^{1+\varepsilon}<\infty$ if
$\varepsilon\in[0,1)$  and $\sup_{n \ge
0}\bar{\E}|X_{n-1}|^{1+\varepsilon}<\infty$ is implied by
$|a|<1$ and $\int_\Real |x|\bar{\nu}(dx)<\infty$. Thus, by Corollary
\ref{theo-2}
$$
\lim_{n\to\infty}\E\Big|\int_\Real |x|
\pi_n(dx)-\int_\Real |x| \bar{\pi}_n(dx)\Big|=0.
$$

\subsection{Additive observation noise}

Suppose
\begin{equation}\label{hhhY}
Y_n = h(X_{n-1}) + \xi_n,
\end{equation}
where $h$ is a fixed measurable function, $\xi=(\xi_n)_{n\ge 1}$
is an i.i.d. sequence of random variables independent of $X$.
Since
$$
\E\big(g(Y_n)|\F^Y_{n-1}\big)=\pi_{n-1}(h)+\E
\xi_1
$$
and if one of the integrability conditions in Theorem \ref{lem-1} is satisfied
for $g(y):=y$, \eqref{wstab} holds true for $h$:
\begin{equation}\label{alaCOC}
\lim_{n\to\infty} \E\big|\pi_n(h)-\bar{\pi}_n(h)\big| =0.
\end{equation}

\begin{rem}
\eqref{alaCOC} resembles the result of J.M.C. Clark et al \cite{COC}
in the continuous time setting: for a general Markov signal
$X=(X_{t})_{t\ge 0}$ and the observations $Y=(Y_t)_{t\ge 0}$ of
the form \eqref{KKCL},
$$
\E\int_0^\infty \big(\pi_t(h)-\bar{\pi}_t(h)\big)^2 dt\le
2\int_\Real \log \frac{d\nu}{d\bar{\nu}}(x)\nu(dx).
$$
This bound is verified information theoretical arguments.
\end{rem}

\subsubsection{\bf Linear observations $\mathbf {h(x)\equiv x}$}
Consider the linear observation model \eqref{hhhY} with $h(x)=x$:
$$
Y_n = X_{n-1} + \xi_n.
$$

\begin{prop}\label{pro-5}
Assume
\begin{enumerate}
\renewcommand{\theenumi}{A\arabic{enumi}}
\item\label{A1}
$\dfrac{d\nu}{d\bar{\nu}}\le c<\infty$
\item \label{A2}
$X_n^p$ is $\bar{\P}$-uniformly integrable for some $p\ge 1$
\item
\label{A3} $\big|\E e^{i\xi_1 t}\big|>0$ for all $t\in \Real$
\end{enumerate}
Then for any continuous function $f(x), x\in\Real$, growing
not faster than a polynomial of order $p$,
$$
\lim_{n\to\infty}\E|\pi_n(f)-\bar{\pi}_n(f)|=0.
$$
\end{prop}

\medskip
\noindent
{\bf Proof.}
If $f$ is an unbounded function, it can be approximated by a sequence of bounded functions
$f_\ell$, $\ell\ge 1$ with $f_\ell(x)=g_\ell(f(x))$, where
$
g_\ell(x)=
  \begin{cases}
    x , & |x|\le\ell \\
    \ell\sign(x) , & |x|>\ell.
  \end{cases}
$

Further, for $k=1,2,\ldots$, set
$$
f_{\ell,k}(x)=
  \begin{cases}
    f_\ell(x), & |x|\le k-1
    \\
    \tilde{f}_{\ell,k}(x), & k-1<|x|\le k
    \\
    0, & |x|>k,
  \end{cases}
$$
where $\tilde{f}_{\ell,k}(x)$ is chosen so that the function
$f_{\ell,k}(x)$ is continuous and
$$
|f_{\ell,k}(x)|\le|f_\ell(x)|.
$$
By the second Weierstrass approximating theorem
(see e.g. \cite{Szego}) one can choose a trigonometrical polynomial
$P_{m,\ell,k}(x)$ such that for any positive number $m$,
$$
\max_{x\in[-k,k]}\big|f_{\ell,k}(x)-P_{m,\ell,k}(x)\big|\le \frac{1}{m}.
$$

Since $P_{m,\ell,k}(x)$ is a periodic function,
$$
|P_{m,\ell,k}(x)|\le \frac{1}{m}+\max_{|y|\le k}|f_{\ell,k}(y)|\le
\frac{1}{m}+\ell, \quad\text{for any $|x|>k$}.
$$
Using the triangular inequality for
\begin{gather*}
f=P_{m,\ell,k}+\big(f_{\ell,k}-P_{m,\ell,k}\big)+\big(f_\ell-f_{\ell,k}\big)+\big(f-f_\ell\big),
\end{gather*}
and the estimates
\begin{list}{}{}
\item $|f_{\ell,k}-P_{m,\ell,k}|\le \frac{1}{m}I_{\{|x|\le k\}}+\Big(\frac{1}{m}+\ell\Big)I_{\{|x|> k\}}$
\item $|f_\ell-f_{\ell,k}|\le \ell I_{\{|x|> k\}}$
\item $|f-f_\ell|\le C(1+|x|^p)I_{\{C(1+|x|^p)> \ell\}}$, for some constant $C>0$
\end{list}

\smallskip
\noindent
we find the following upper bound
\begin{gather*}
|f-P_{m,\ell,k}|\le \frac{1}{m} I_{\{|x|\le k\}}+ \Big(\frac{1}{m}+2\ell\Big)I_{\{|x|> k\}}+
C(1+|x|^p)I_{\{C(1+|x|^p)> \ell\}},
\end{gather*}
implying
\begin{multline*}
\E\big|\pi_n(f)-\bar{\pi}_n(f)\big| \le
\E\big|\pi_n(P_{m,\ell,k})-\bar{\pi}_n(P_{m,\ell,k})\big|
+\frac{2}{m}+2\ell\E\int_{\{|x|>k\}}[\pi_n(dx)+\bar{\pi}_n(dx)]
\\
+ C\E\int_{\{C(1+|x|^p)\ge \ell\}} (1+|x|^p)\pi_n(dx)
+C\E\int_{\{C(1+|x|^p)\ge \ell\}} (1+|x|^p)\bar{\pi}_n(dx).
\end{multline*}
Thus, the desired result holds by arbitrariness of $m$ provided that

\begin{enumerate}
\renewcommand{\theenumi}{(\arabic{enumi})}
\item \label{raz} $\lim\limits_{n\to\infty}\E\big|\pi_n(P_{m,\ell,k})-\bar{\pi}_n(P_{m,\ell,k})\big|=0$,
$\forall \ m,\ell,k $;

\item \label{dva}  $\lim\limits_{k\to\infty}\varlimsup\limits_{n\to\infty} 2\ell\E\int_{\{|x|>k\}}
[\pi_n(dx)+\bar{\pi}_n(dx)]=0$, $\forall \ \ell$;

\item \label{tri} $\lim\limits_{\ell\to\infty}\varlimsup\limits_{n\to\infty}
\E\int_{\{C(1+|x|^p)\ge \ell\}} (1+|x|^p)[\pi_n(dx)+\bar{\pi}_n(dx)]=0$;
\end{enumerate}

\medskip
\noindent
\ref{raz} holds due to
$
\E\big(e^{i tY_n}\big|\F^Y_{n-1}\big)=
\E\big(e^{iX_{n-1} t}\big|\F^Y_{n-1}\big)\E e^{it\xi_1}
= \pi_{n-1}(e^{itx})\E e^{it\xi_1}
$
and the assumption \eqref{A3} since,  by Theorem \ref{lem-1},
$$
\lim_{n\to\infty}\E\big|\pi_n\big(e^{itx}\big)-\bar{\pi}_n\big(e^{itx}\big)\big|=0,
\quad \forall t\in\Real.
$$

\medskip
\noindent
\ref{dva} is implied by the Chebyshev inequality
$$
\E\int_{\{|x|>k\}}
[\pi_n(dx)+\bar{\pi}_n(dx)]\le\frac{1}{k}\bar{\E}
\Big(1+\frac{d\nu}{d\bar{\nu}}(X_0)
\Big)|X_n|,
$$
and the assumptions \eqref{A1} and \eqref{A2}.

\medskip
\noindent
\ref{tri} follows from
$$
\E\int_{\{C(1+|x|^p)\ge \ell\}} (1+|x|^p)[\pi_n(dx)+\bar{\pi}_n(dx)]=
\bar{\E} I_{\{C(1+|X_n|^p)\ge\ell\}}\Big(1+\frac{d\nu}{d\bar{\nu}}(X_0)
\Big)\big(1+|X_n|^p)
$$
and the assumptions  \eqref{A1} and \eqref{A2}.\hfill $\blacksquare$

\subsection*{Acknowledgement.} We are grateful to Valentine Genon-Catalot
for bringing \cite{G-C}, \cite{G-C2} to our attention.

\end{document}